\documentclass[12pt]{article}

\usepackage{amsmath}
\usepackage{amssymb}
\usepackage{epsfig}
\usepackage{psfrag}
\usepackage{color}

\usepackage{geometry}
\newcommand{\R}{ {\mathbb R} }

\newcommand{\T}{{\mathbb T}}

\DeclareMathOperator{\supp}{Supp}

\newcommand{\comments}[1]{}
\newtheorem{theorem}{Theorem}[section]

\newtheorem{prop}[theorem]{Proposition}
\newtheorem{lemma}[theorem]{Lemma}

\newtheorem{remark}[theorem]{Remark}

\makeatletter
\@addtoreset{equation}{section}
\makeatother

\newcommand{\cqfd}{{\unskip\kern 6pt\penalty 500
\raise -2pt\hbox{\vrule\vbox to 6pt{\hrule width 6pt
\vfill\hrule}\vrule}\par}}

\begin{document}
\title{Stability of trajectories for $N$-particles dynamics with 
singular potential.}
\author{J.~Barr\'e$^1$,\quad jbarre@unice.fr\\
M.~Hauray$^2$,\quad hauray@cmi.univ-mrs.fr\\ 
P.~E.~Jabin$^{1,3}$,\quad jabin@unice.fr} 
\footnotetext[1]{Laboratoire J.~A.~Dieudonn\'e, UMR CNRS 6621,
Universit\'e de Nice-Sophia Antipolis, Parc Valrose, 06108 Nice France. }
\footnotetext[2]{Centre de Math\'ematiques et Informatique (CMI),
  Universit\'e de Provence, 
Technop\^ole Ch\^ateau-Gombert, 39, rue F. Joliot Curie, 13453 Marseille
Cedex.}
\footnotetext[3]{TOSCA project-team, INRIA Sophia Antipolis --
  M\'editerran\'ee, 2004 route des Lucioles, B.P.\ 93, 06902 Sophia
  Antipolis Cedex, France.}
\date{}
\maketitle
\begin{abstract}
We study the stability in finite times of the trajectories of interacting
particles. Our aim is to show that in average and {\em uniformly in
  the number of particles}, two trajectories whose
initial positions in phase space are close, remain close enough at
later times. For potential less singular than the classical
electrostatic kernel, we are able to prove
such a result,  for initial positions/velocities distributed
according to the Gibbs equilibrium of the system. 
\end{abstract}

%%%%%%%%%%%%%%%%%%%%%%%%%%%%%%%%%%%%%%%%%%%%%%%%%%%%%%%%%%%%%%%%%%%%%%%
\section{Introduction}
The stability of solutions to a differential system of the type
\begin{equation}
\frac{dZ}{dt}=F(Z(t)),\label{odegeneral}
\end{equation}
is an obvious and important question. For times of order $1$ and if
$F$ is regular enough (for instance uniformly Lipschitz), 
the answer is given quite simply by Gronwall lemma. For two
solutions $Z$ and $Z^\delta$ to \eqref{odegeneral}, one has 

\begin{equation}
|Z(t)-Z^\delta(t)|\leq |Z(0)-Z^\delta(0)|\,\exp(t\,\|\nabla F\|_{L^\infty}).
\label{gronwall}
\end{equation}
This inequality forms the basis of the classical Cauchy-Lipschitz
theory for the well posedness of \eqref{odegeneral}. 
It does not depend on the dimension of the system (the norm chosen is
then of course crucial). This is hence very convenient for the study
of systems of interacting particles, which is our purpose here.

Consider the  system of equations
\begin{equation}
\left\{ \begin{array}{l}
\dot{X}^N_i=V^N_i\\
\dot{V}^N_i= E_N(X^N_i) = \frac{1}{N}\sum_j K(X^N_i-X^N_j)
\end{array} \right.
\label{eqs:npart}
\end{equation}
where for simplicity
all positions $X^N_i$ belong to the torus $\T^3$ and all
velocities $V^N_i$ belong to $\R^3$. This system is obviously 
a particular case
of \eqref{odegeneral} with $Z=Z^N=(X_1^N,\ldots,X_N^N,V_1^N,\ldots, V_N^N)$.  

The equivalent of \eqref{gronwall} reads in this case
\begin{equation}
\|Z^N(t)-Z^{N,\delta}(t)\|_1\leq \|Z^N(0)-Z^{N,\delta}(0)\|_1\,
\exp(t\,(1+\|\nabla K\|_{L^\infty})),\label{gronwallN}
\end{equation}
where we define
the norm on $\Pi^{3N}\times \R^{3N}$
\[
\|Z\|_1=\frac{1}{2N}\,\sum_{i=1}^N (|X_i|+|V_i|).
\]
This estimate is logically uniform in the number of particles $N$. It
is important in itself but also because it is a crucial tool to pass
to the limit in the system of $N$ particles and derive the Vlasov-type
equation 
\begin{equation}
\partial_t f+v\cdot\nabla_x f+\left(F\star_x\left(\int_{\R^3}
f(t,x,v)\,dv\right)\right)\;\cdot \nabla_v f=0,\label{vlasov}
\end{equation}
for the 1-particle density $f(t,x,v)$ in phase space, where 
$\star$ denotes the convolution.
Hence estimates such as \eqref{gronwallN} are at the heart of the
derivation performed in \cite{Neunzert},
\cite{Dobrushin}, and \cite{BH77} (we also
refer to \cite{Sp}, \cite{VA} and \cite{Wo}). Note that the derivation
of collisional kinetic
models (of Boltzmann type) involves quite different techniques, see
\cite{Lanford} or \cite{CIP}. 

Unfortunately, many cases of interest in physics deal with singular
forces $K\not\in W^{1,\infty}_{loc}$. Typical cases are $K= -\nabla
\phi$, a periodic force coming from a periodisation of the potential
$\phi_{\R}(x)\sim C/|x|^{\alpha-1}$, {\em i.e.} \mbox{$ \phi(x) =
  \phi_{\R}(x) + g(x)$,} where $g(x)$ is an (at least) $C^2$ function
on the torus $\T^3$. As the potential $\phi$ is defined up to a
constant, we may also assume that its average is $0$: $\int_{\T^3}
\phi(x)\,dx =0$. The most important case is the electrostatic or
gravitational interaction: $\alpha=2$, in dimension $3$.

Very little is known for these singular kernels, either from the point
of view of the stability or of the derivation of Vlasov-type equations. 
Provided $\alpha<1$ and
the initial configuration of particles are well distributed, the limit
to Vlasov equation \eqref{vlasov} was proved in \cite{HJ06}. 

For systems without inertia, {\em i.e.} when the equations are simply
\begin{equation}
\left\{ \begin{array}{l}
\dot{X}^N_i= E_N(X^N_i) = \frac{1}{N}\sum_j K(X^N_i-X^N_j)
\end{array} \right.,
\label{eqs:npartmacro}
\end{equation}
it seems to be easier to implement Gronwall-type inequalities. The
derivation of the mean field limit is consequently known up to
$\alpha<2$ ($\alpha<d-1$ in dimension $d$), see \cite{Hau2} and also
\cite{JO} for a situation where the forces have a
more complicated structure. In this
setting the most important case is however found in dimension 2, for
$K=x^\perp/|x|^2$ (corresponding to $\alpha=d-1=1$);
the limit is the 2d incompressible Euler
equation written in vorticity form. The derivation of the mean-field
limit in this case  
was rigorously performed in \cite{GHL90} and \cite{Sch1}, \cite{Sch2}.

For differential equations like \eqref{odegeneral} in finite
dimensions, it has long been known that well posedness and stability
(for almost all initial data)
can be achieved without using Gronwall-type estimates. The
introduction of renormalized solutions by DiPerna-Lions in \cite{DL}
gave well posedness for $F\in W^{1,1}$ with $div\,F\in L^\infty$. 

This was extended to $F\in BV$ in the phase space situation in
\cite{Bo} and then in the general case in \cite{Am} (see also
\cite{HLL} for a slightly different approach). The exact case of the
Poisson interaction was treated in \cite{Hau}.

This well posedness
implies some stability as the flow has then some differentiability
properties, see \cite{ALM}. However the corresponding stability
estimate is not quantitative and this kind of method does not seem to
be able to provide uniform estimates
in the number of particles (which gives the dimension of the
system). We refer to \cite{DeL} for a precise overall presentation of
the well posedness and differentiability issues for
Eq. \eqref{odegeneral} in finite dimension.

More recently a new method to show well posedness for
\eqref{odegeneral} has been introduced in \cite{CD}. Given a fixed
shift $\delta$, it consists in bounding quantities like
\begin{equation}
\int_{Z^0}
\log\left(1+\frac{|Z(t,Z^0)-Z(t,Z^0+\delta)|}{|\delta|}\right)\,dZ^0, 
\label{log}
\end{equation}
where $Z$ is the flow associated to \eqref{odegeneral}, {\em i.e.}
$Z(t,Z^0)$ is a solution to \eqref{odegeneral} satisfying $Z(0,Z^0)=Z^0$. 

A bound on such a quantity shows that for {\em a.e.} $Z^0$ the two
trajectories $Z(t,Z^0)$ and $Z(t,Z^0+\delta)$ remain at a distance of
order $|\delta|=|Z(0,Z^0)-Z(0,Z^0+\delta)|$. In this sense, this is an
almost-everywhere version of the Gronwall inequality \eqref{gronwall}.

It was shown in \cite{CD} that the quantity \eqref{log} remains
bounded if $F\in W^{1,p}$ for some $p>1$. This was extended to
$W^{1,1}$ and $SBV$ in \cite{Ja} and even to $H^{3/4}$ in the phase
space setting in \cite{CJ}. However in all those results the bound
depends on the dimension of the space and is blowing-up as this
dimension increases to $\infty$.

Therefore, 
our precise aim in this article is to prove a bound on a quantity like
\eqref{log} for the system \eqref{eqs:npart}, uniformly in the number of
particles. To our knowledge, this is the first quantitative stability
estimate to be obtained for singular forces.

Several new important issues occur when one tries to do that
though. One of the most important is the reference measure which is
chosen as this can imply different notions of almost everywhere
as the dimension tends to $+\infty$.  
In finite dimension, this refers to the Lebesgue
measure and of course implies corresponding estimates for any measure
which is absolutely continuous with respect to the Lebesgue
measure. In infinite dimension, no such natural measure exists. 
This is due to
the phenomenon of concentration of measures. Even in the case of
finite but large dimensions, this is a problem to get quantitative
estimates. Indeed even if two measures $\nu_1$ and $\nu_2$ on
$\Pi^{3N}\times \R^{3N}$ are absolutely continuous with respect to
each other or even more if
\[
d\nu_1\leq C\,d\nu_2,
\]  
then the constant $C$ will in general depend on the dimension and 
go to $+\infty$ as $N$ increases. This means that a uniform
quantitative estimate on the trajectories for some measure $d\nu_1(Z^0)$
on the initial configuration will not in general imply a good estimate for
another measure $d\nu_2(Z^0)$.

\medskip

For each $N$, the choice of the measure $\mu_N$ on $\Pi^{3N}\times
\R^{3N}$ is hence crucial to get a good estimate. One would naturally
want to select a measure $\mu_N$ which is bounded, stable and invariant by the
flow, just as the Lebesgue measure is stable and invariant by the flow of
\eqref{odegeneral} when $F$ is divergence free. 
Let us denote by $Z^N(t)= (X^N(t),V^N(t))$ or $Z^N(t,Z_0^N)$ 
the vector of particles velocities
and positions evolving through Eq. \eqref{eqs:npart} till time $t$ and
depending on the initial configuration $Z_0^N$. 
The system \eqref{eqs:npart} has an invariant which is the total energy
\begin{eqnarray}
H_N[Z^N] & = & \sum_{i}\frac{(V^N_i)^2}{2}+\frac{1}{2N}\sum_{i\neq
  j}\phi(|X^N_i-X^N_j|)\\  % -(N-1)e_0~, 
& = & E_{kin}(V^N) + E_{pot}(X^N)
\end{eqnarray}
To get an invariant measure $\mu_N$, the simplest choice is to take a
function of the total energy. Among those which are stable, the most
natural is the Gibbs equilibrium
\begin{equation}
d\mu_N(Z^N)=\frac{1}{\mathcal{B}_N}e^{-\beta H_N[Z^N]} dZ^N~,
\end{equation}
where $dZ^N$ is Lebesgue measure on $\T^{3N} \times \R^{3N}$, and
\begin{equation} \label{eq:betaN}
\mathcal{B}_N(\beta)=\int e^{-\beta H_N[Z^N]}~dZ^N~,
\end{equation}
is the normalization constant. Note that for a potential $\phi$ with a
singularity at the origin, this makes sense only if $\lim_0 \phi =
+\infty$, that is in the case of repulsive interactions.  In the
following, since we deal with measures which have a density with
respect to the Lebesgue measure $dZ_0^N$, we will use the same
notation for the measure $\mu_N$ and its density.  

\medskip

We study the quantity
\begin{equation} \label{eq:defQ}
Q(t)=\int~d\mu_N(Z^N_0)\int_{\delta\in \T^{3N}\times
  \R^{3N}}\;\psi_N(Z^N_0,\delta)
\ln \left(1+\frac{\|Z^N(t,Z^N_0)-Z^N(t,Z^N_0+\delta)\|_1}{\delta_N}
\right)\,d\delta ~,
\end{equation}
where $\delta_N$ is a small parameter that will go
slowly to zero when $N$ goes to infinity.
$\delta$ is a shift on the initial condition $Z^N_0$

% and which gives the order of
% magnitude of the shift $\delta$.

Here $\psi_N : \T^{3N} \times \R^{3N} \mapsto \mathcal{P}(\T^{3N} \times
\R^{3N})$ (where $\mathcal{P}(\Omega)$ denotes the set of probabilities 
on $\Omega$) is a probability valued function, so that it satisfies
\[
\int_{\delta \in \T^{3N} \times \R^{3N}}  d \psi_{N}(Z_0^N,\delta)=1, \quad
\forall Z_0^N \in (\T^3 \times \R^3)^N
\]
$\psi_N$ then then gives distribution of the shifts on the initial 
conditions, and the quantity $Q(t)$ is averaged both on the initial 
conditions $Z^N_0$ and on the shifts $\delta$.

We now define the image measure of $\mu_N$ by the shift distribution:
\[
\tilde\mu_N(Z_0) =  \int \mu_N(Z_0-\delta) \psi_N(Z_0-\delta ,\delta)
\,d\delta \,.
\] 
The crucial assumption will be that the image measure $\tilde{\mu}_N$ 
remains ``close'' to the original Gibbs measure in the sense that
\begin{equation} \label{psi}
\exists \,K_\beta > 0, \text{ such that } \forall\, Z_0 \text{ and } N, 
\quad\tilde{\mu}_N(Z_0) \leq K_\beta \mu_N(Z_0) 
\end{equation}
with a constant $K_\beta$ independent of $N$, but which may depend on $\beta$. 

We will also use the weaker but very similar condition
\begin{equation} \label{psi2}
\exists \,K'_\beta >0 \text{ s.t. } \forall \, N, \exists \, \beta'
\leq \beta, \text{ s.t. } \forall\, Z_0,
\;\tilde{\mu}_N(Z_0) \leq K'_\beta \mu_N^{\beta'}(Z_0),
\end{equation}
where $\mu_N^{\beta'}$ denotes the Gibbs measure with inverse temperature
$\beta$. The last condition is more general than the first, it allows to
control the image measure by a Gibbs measure with bigger temperature. 

\begin{remark}
 We mention that by definition of $\tilde\mu_N$, the measure $\pi_N(Z_0,Z_0') =
\mu_N(Z_0) \Phi_N(Z_0, Z_0' -Z_0)$ is a transport from the measure $\mu_N$ to
$\tilde\mu_N$. In fact, $\psi_N$ is (up to a translation of origin) what is
usually called the desintegration of the measure $\pi_N$ with respect to its
first projection $\mu_N$. However, we preferred our less standard presentation
since we are more interested in $\mu_N$ and its shift $\psi_N$ that in the
precise image measure $\tilde\mu_N$. We mention the analogy to
emphasize that our quantity $Q$ is related to optimal transport. In fact
\[
 Q_N(0) = \int~d\mu_N(Z^N_0)\int_{\delta\in \T^{3N}\times
  \R^{3N}}\;d\psi_N(Z^N_0,\delta)
\ln \left(1+\frac{\|\delta\|_1}{\delta_N}
\right) \geq W_N(\mu_N,\tilde\mu_N)
\]
where $ W_N$ is the transport for the cost $\ln\left( 1 +
\frac{\|\cdot\|_1}{\delta_N}\right)$.
\end{remark}

Conditions \eqref{psi} and \eqref{psi2} are not explicit on $\psi$.
Roughly speaking, they will be satisfied if $\psi_N$ is chosen so
that $|H_N(Z_0+\delta) - H_N(Z_0)| \leq C$ if $\delta \in \supp
\Psi_N(Z_0,\cdot)$.  This is reasonable since $H_N$ is preserved by
the dynamics, so that if the shift changes the energy too much, the
original and shifted dynamics may be very different. But that simple
and ``reasonable'' condition is not sufficient, we will really need a
bound like \eqref{psi} on the image measure constructed with the
shift.

As the conditions are not explicit, we will provide in section \ref{examples}
some examples of admissible shift distributions. 
The main result of the paper is a control on the growth on this
quantity $Q$: 
%%%%%%%%%%%%%%%%%%%%%%%%%%%%%%%%%%%%%%%%% THM1
\begin{theorem} \label{theorem1}
Assume that $\phi  \geq \phi_{min}$ for some $\phi_{min} \in \R^-$ and
that for some constant $C$, and  $\alpha<2$ 
\[
 \phi(x)\leq \frac{C}{|x|^{\alpha-1}},\quad |\nabla \phi|\leq
\frac{C}{|x|^{\alpha}},\quad  |\nabla^2 \phi|\leq\frac{C}{|x|^{\alpha+1}}.
\]
Then taking $\delta_N=N^{-\varepsilon}$ for any
$\varepsilon \leq 1-\alpha/3$ and for all $N \geq \frac{6^4}{(2-\alpha)^2}$
one has
\[
Q(t)\leq   \left(1+ (1+ K_\beta)\frac{C c_\beta  +C_a
  c_\beta^a}{2-\alpha} \right)\,t + Q(0) ~,  
\]
where 
%$c_\beta = (\beta^{3/2}e^{\beta e_0})/(2\pi)^{3/2}$, 
$c_\beta = e^{- \frac \beta 2 \phi_{min}}$, 
$a$ is any exponent strictly larger than $2\alpha/3$, $C$ constant
(that can be made explicit),  
and $C_a$ satisfies $C_a \leq \frac{C}{3a - 2\alpha}$.
\end{theorem}

This theorem is not able to deal with the electrostatic interaction,
$\alpha=2$; gravitational is of course out of question since repulsive
potentials are needed. Note however that the electrostatic potential
is just the critical case. The same result could be obtained in any
dimension, with essentially the same proof. In dimension $d$, the
condition would then be $\alpha<d-1$. The growth of $Q$ is linear in
time: note that this indeed corresponds to an average exponential in
time divergence of the trajectories, analogous to \eqref{gronwall}.

\medskip

Roughly speaking, and provided that the average shift at time $0$ is
of order $\delta_N$ (or smaller), the theorem says that the average
shift transported by the dynamics remains of order $\delta_N$ during
the evolution, and the control given is quite good.  It is interesting
to compare $\delta_N$ to the minimal distance in the $(X,V)$ space
between two particles of a configuration, which is of order
$N^{-\frac13}$. Notice that it is always possible to choose $\delta_N$
smaller than $N^{-\frac13}$. Then if the order of magnitude of the
initial shifts is smaller than $N^{-\frac13}$, the theorem says it
remains so at all time. This implies that the pairing of a particle in
the configuration $Z(t)$ with the closest one in the configuration
$Z^{\delta}(t)$ is not very much affected by the dynamics: in this
sense, there is not much ``mixing''.

%We may compare it
%to the minimal distance between the particles of a configuration which
%is of order  
%$N^{\frac13}$ in $(x,v)$ space. So when we choose a configuration
%$Z_0$ and shift it to a configuration $Z_0^\delta$ with a shift of
%order $N^{\frac12}$, we may say that the shift associate to a particle
%$Z_0$ to its closest particle in the configuration $\delta$. Since
%$\delta_N$  may always be chosen smaller than $ N^{\frac13} $ our
%theorem says that it remains true during all time. At time $t$, the
%shift still associate to most of the particle of $Z(t)$ its closest
%one in $Z^\delta(t)$.  So that there is not much mixing.  

\medskip

While the Gibbs equilibrium is the most natural choice for the measure
$\mu_N$, others are possible. The proof would work for any measure
$\mu_N$ such that 
\begin{itemize}
\item[-] $\mu_N$  is invariant under the flow or
  $\mu_N(Z^N(t))=\mu_N(Z_0^N)$ for  $Z^N$ solution to
  \eqref{eqs:npart} 
\item[-] for all $k$, the  $k$-marginal of $\mu_N$ defined by $
  \mu_N^k(Z^k)=\int \mu_N(Z_0^k,\tilde{Z_0}^{N-k}) d\tilde{Z_0}^{N-k}
  $  satisfies:   
\[
\forall Z_0, \quad \mu_N^k(Z_0) \leq C^k.
 \]
\end{itemize}

 Obvious candidates are functions of the renormalized energy $\frac{1}{N}H_N$
but checking the bound on the marginals is not necessarily easy. 

\medskip
{ \bf  Link with mean field limit. }

Finally, let us emphasize that this stability estimate does not answer
the question of the mean field limit. Doing so would require to be
able to deal with much more general measures $\mu_N$. More precisely if
one can prove Th. \ref{theorem1} for a sequence of $\mu_N$ and if in
some reasonable sense
\[
\mu_N-\Pi_{i=1}^N f^0 \longrightarrow 0,
\]
then the mean field limit is proved but only for the initial data
$f^0$. Currently the Gibbs equilibrium corresponds to
$f^0(x,v)=e^{-\beta |v|^2/2}$.

 Unfortunately, we are not able do deals with more general measures
$\mu_N^0$. The problem is that we need bounds on every $k$
marginals and those are very difficult to obtain if we start from another
measure than the Gibbs equilibrium. For instance, starting from $\mu_N =
g^{\otimes N}$ for some smooth profil $g$, we have the desired bounds
at time $0$, but do not know if $\mu_N(t)$ satisfies them for any other
time $t> 0$.
   
%{ \bf  M : Je souhaiterais qu'on soit plus precis dans ce passage. Par
%  exemple. P-E: Je prefere le debut du passage en haut. Il me semble
%  que tu veux surtout parler de marginales. Donc je te propose de
%  rajouter quelques lignes la-dessus, en italique plus haut. }

%Finally, let us emphasize that this stability estimate does not answer
%the question of the mean field limit. Doing this will means to obtain
%a result of the following type. If a measure $\mu_N^0$ on the $N$
%particles configuration space concentrate on a smooth profile $f$ on
%$\T^3 \times \R^3$ (in the sense that when we choose the empirical
%measures $\frac 1 N \sum \delta_{(X_i,V_i)}$ according to $\mu_N$,
%they are close to $f$ in some sense), does the measure $\mu_N(t)$
%concentrate on the profile $f(t)$.  

%When $\mu_N$ is a Gibbs measure, such a result is trivial, since  this
%measure concentrate on the profil $f^0(x,v)=e^{-\beta |v|^2/2}$  (See
%ref???). But at time $t$, $\mu_N(t) =\mu_N$ and $f(t)=f^0$ since both
%are invariants, so that there is nothing to prove. 

%And unfortunately, we are not able do deals with more general measures
%$\mu_N^0$. The problem is that we need good bounds on the $k$
%marginals, that are very difficult to prove if we start form another
%measure than the Gibbs one. For instance, starting from $\mu_N =
%g^{\otimes N}$ for some smooth profil $g$, we have the desired bound
%at time $0$, but do not know if $\mu_N(t)$ satisfy them for another
%time $t \neq 0$.   

\section{Some examples of admissible shift distributions.} 
\label{examples}

\subsection{Shift on velocity variables} \label{speedshift}

In this section, we will be interested  in shifts acting only the
velocities. A first possibility is to take shifts independent of $Z_0$ and
acting independently and identically on each velocity
variable. Precisely, we are looking for shift distributions such as 
\[
\psi_N(\delta)  =  \delta_{\delta_X =0} \prod_{i=1}^N N^{\frac32} \psi(\sqrt N
\delta_{v_i} )                  
\]
where $\psi$ is a probability on $\R^3$ symmetric with respect to the origin
($\psi(-v) = \psi(v)$ if $\psi$ has a density). 

We will not be able to deal with a general $\psi$, but will show that the
hypothesis \eqref{psi2} is satisfied if $\psi$ is Gaussian or has a compact
support. This is stated precisely in the following Proposition:

\begin{prop} \label{prop:spshift}
Assume that $\psi$ is a Gaussian probability with variance $\sigma^2$ :
\[ \psi(\delta_v)=\left(\frac{1}{2\pi\sigma^2}\right)^{\frac32}
e^{-\frac{\delta_v^2}{2\sigma^2}}
\]
then the hypothesis \eqref{psi2} is satisfied with 
\[
 \beta'(N) = \beta\left(1-\frac{1}{1+N/(\beta\sigma^2)}\right) \,, 
\text{ and }
 K_\beta = e^{-\frac{\beta^2\sigma^2}{4}\phi_{min}}
e^{\frac34 \beta\sigma^2}\,.
\]
Assume that $\psi$ has a compact support with $\supp \psi \subset B(0,\delta_m)$
(the ball of radius $\delta_m$, center $0$). Then \eqref{psi2} is satisfied for
$N > \beta \delta_m^2$ with 
\[
 \beta'(N) = \beta \left( 1  - \frac{\beta \delta_m^2}{N}\right) \,, 
\text{ and }  K_\beta = e^{-
  \frac{\beta^2 \delta_m^2 \phi_{min}}{2} } e^{\frac32 \beta\sigma^2}
\]

\end{prop}

Remark that for $\alpha \leq \frac32$, the average velocity fluctuation
given by such shift is larger than the smallest $\delta_N$ we can
choose.

\medskip

{\bf  Proof of the proposition \ref{prop:spshift}}\\
In the Gaussian case, we have
\[
\tilde{\mu}_N(Z_0) = \frac{e^{- \beta H_N(Z_0)}}{\mathcal{B}_N(\beta)}
 \left( \frac{N }{2\pi\sigma^2}\right)^{\frac{3N}2}
\int e^{- \frac\beta{2} \sum_i (\delta_{v_i}^2 - 2 v_i
  \delta_{v_i})}  e^{-\frac{N}{2\sigma^2}\sum_i\delta_{v_i}^2 } \,d\delta_V  \,.
\]
The $3N$ integrals may be performed independently, using the 1D
calculation:
\[
 \sqrt{\frac{N}{2\pi\sigma^2}} \int_{-\infty}^{+\infty}
 e^{-\frac12 (\beta+N/\sigma^2)\delta_v^2 +\beta v\delta_v} \,d\delta_v
= \sqrt{\frac{N}{N+\beta\sigma^2}}
e^{\frac\beta2 \frac{\beta\sigma^2}{N+\beta\sigma^2}}
\]
We finally get
%\[
\begin{eqnarray}
\tilde{\mu}_N(Z_0) &=& \frac{e^{-\beta H_N(Z_0)}}{\mathcal{B}_N(\beta)}
 \left( \frac{1 }{1+\beta\sigma^2/N}\right)^{\frac{3N}2}
e^{\frac{\beta}{1+N/\beta\sigma^2}E_{kin}(Z_0)} \nonumber \\
&=&  \frac{e^{-\beta' H_N(Z_0)}}{\mathcal{B}_N(\beta')}
 \frac{\mathcal{B}_N(\beta')}{\mathcal{B}_N(\beta)}
\left( \frac{1 }{1+\beta\sigma^2/N}\right)^{\frac{3N}2}
e^{-\frac{\beta}{1+N/\beta\sigma^2}E_{pot}(Z_0)} \nonumber  
\end{eqnarray}
%\]
with $\beta'(N)=\beta(1-\frac{1}{1+N/(\beta\sigma^2)})$. Using
\[
\frac{\beta}{1+N/\beta\sigma^2}E_{pot}(Z_0)\geq 
\frac{\beta^2\sigma^2\phi_{min}}{4}~,
\]
we get
\[
\tilde{\mu}_N(Z_0)\leq e^{-\frac{\beta^2\sigma^2}{4}\phi_{min}}
e^{-\frac34 \beta\sigma^2} \frac{\mathcal{B}_N(\beta')}{\mathcal{B}_N(\beta)}
\tilde{\mu}_N^{\beta'}(Z_0)
\]
We define $\mathcal{B}_{N,X}$ the normalization integral restricted to the 
position variables:
\[
\mathcal{B}_{N,X}(\beta)= \int e^{-\beta E_{pot}(X_1,\ldots,X_N)}\,dX_1\ldots dX_N
\]
From the Jensen inequality applied to the function 
$x\mapsto x^{\beta'/\beta}$, we obtain for $\beta' \leq \beta$  
\[
\left( \mathcal{B}_{N,X}(\beta')\right)\leq 
\left( \mathcal{B}_{N,X}(\beta)\right)^{\beta'/\beta} 
%\left( \frac{2\pi}{\beta'}\right)^{3N/2}
%\left( \frac{\beta_2}{2\pi}\right)^{(3N/2)(\beta_1/\beta_2)}
\]
Then, using the bounds of Lemma \ref{lemma2}, 
we have $\mathcal{B}_{N,X}(\beta)^{\beta'/\beta-1}\leq 1$; this implies
\[
 \frac{\mathcal{B}_N(\beta')}{\mathcal{B}_N(\beta)}\leq 
\left( \frac{\beta}{\beta'}\right)^{3N/2}
\leq e^{3\beta\sigma^2/2} 
\]
This proves the desired inequality with
\[
 K_\beta = e^{-\frac{\beta^2\sigma^2}{4}\phi_{min}}
e^{\frac34 \beta\sigma^2}\,.
\]
\medskip

In the case of $\psi$ with compact support, we follow the same sketch.
To do this, we will need a bound on
\[
 \int e^{\frac \beta 2 (2 v \delta_v -
\delta_v^2)}N^{\frac32}\,d\psi(\sqrt N \delta_v) 
= \int e^{\frac \beta 2 \left( \frac{2 v \delta_v}{\sqrt N} -
\frac{\delta_v^2}N \right)} \,d\psi(\delta_v)  \,.
\]
Using the symmetry of $\psi$, and the inequality $\cosh(x) \leq
e^{\frac{x^2}2}$ (valid for $x \in \R$), we may bound that term by
\[
\int \cosh \left(\frac{\beta  v \delta_v}{\sqrt N} \right)
e^{- \frac{\beta\delta_v^2}{2N}} \,d\psi(\delta_v) 
\leq 
\int 
e^{ \frac{\beta^2  v^2 \delta_v^2}{2N} - \frac{\beta\delta_v^2}{2N}}
\,d\psi(\delta_v)  \leq
e^{ \frac{\beta^2  \delta_m^2 v^2 }{2N}}
\]
and as in the previous case, we get
\[
\int \mu_N(Z_0-\delta)\psi_N(Z_0-\delta,\delta) \,d\delta \leq
\frac{1}{\mathcal{B}_N} e^{-\beta\left( H_N(Z_0) -\frac{\beta \delta_m^2}N
  E_{kin}(Z_0)\right) } \,.
\]
From this point, following exactly the same step as in the case of
gaussian $\psi$, we prove that the hypothesis \eqref{psi2} is
satisfied with the announced constant \cqfd.

\bigskip

By making the shift depend on the velocity $V_0=(V_1,\ldots V_N)$, one
may essentially remove the condition on the size of the norm of the
shift.
More precisely we limit ourselves to shifts $\delta=(0,\ldots,0,\delta')$ with
$\delta'\in\R^{3N}$,
 giving
$Z_0+\delta=(X_0,V_0+\delta')$ with $X_0=(X_1,\ldots, X_N)$. Now define
\begin{equation}
\psi_N(Z_0,\delta)=\left(\Pi_{i=1}^N \delta_{\delta_i=0}\right)\;
 \Psi(|\delta'|)\,G(|V_0|)\,
\delta_{2\,V_0\cdot\delta'+|\delta'|^2=0},\label{shiftv}
\end{equation}
where $|\xi|^2=\xi\cdot\xi$ is the usual euclidian norm and
$\delta_{...}$ is the corresponding Dirac mass on the 
hypersurface of equation
$2\,V_0\cdot\delta'+|\delta'|^2=|V_0+\delta'|^2-|V_0|^2=0$  which is
precisely the sphere of radius $|V_0|$.

We need the
function $\psi_N$ to
satisfy
\[
1=\int
\psi_N(Z_0,\delta)\,d\delta=G(|V_0|)\,\int_{2\,V_0\cdot\delta'+|\delta'|^2=0}
\Psi(|\delta'|)\,d\delta',
\]  
this is always possible with the right choice of $G$ as the integral
\[
\int_{2\,V_0\cdot\delta'+|\delta'|^2=0}
\Psi(|\delta'|)\,d\delta'
\] 
depends only on $|V_0|$ by the rotational symmetry of the sphere.  

Condition \eqref{psi} is automatically true since, as $\mu_N$
depends only on $X_0$ and $|V_0|$ and $|V_0-\delta'|=|V_0|$, 
\[
\tilde\mu_N(Z_0)=\int \mu_N(Z_0-\delta)\,\psi_N(Z_0-\delta,\delta)  
=\mu_N(Z_0)\,\int \psi_N(Z_0,-\delta)=\mu_N(Z_0),\]
because $\psi_N(Z_0-\delta,\delta)=\psi_N(Z_0,-\delta)$.
One could wish to impose additionally that $|\delta'|_1\leq \delta_N$, so 
that $Q(0)$ is of order~$1$. Since
\[
|\delta'|_1=\frac{1}{N}\sum_i |\delta_i'|\leq \frac{1}{\sqrt{N}}\,|\delta'|,
\]
it is enough to impose that $\Psi$ has support in $[0,  N^{1/2}\,\delta_N]$.

As a conclusion, we proved
\begin{prop} For any measure $\Psi$ on $\R^{3N}$, there exists a
  function $G(|V_0|)$, s.t. the probability density $\psi_N$ defined
  by \eqref{shiftv} satisfies \eqref{psi}.
\end{prop}

\bigskip

One could try to generalize this example by letting $V_0+\delta'$ and
$V_0$ to be on close but different energy spheres. For example by posing
\[
\psi_N(Z_0,\delta)=\left(\Pi_{i=1}^N \delta_{\delta_i=0}\right)\;
 \Psi(|\delta'|)\,G_{\eta}(|V_0|)\,
\,{\mathbb I}_{|2V_0\cdot\delta+\delta^2|\leq
  \eta}.
\]
Provided that $\eta$ in not too large, the computations are
essentially the same and one has essentially to make sure that
\[
G_{\eta}(|V_0|\pm \eta)\leq C\, G_{\eta}(|V_0|).
\]  
%%%%%%%%%%%%%%%%%%%%%%%%%%%%%%%%%%%%%%%%
\subsection{Shifts in position variables}
%%%%%%%%%%%%%%%%%%%%%%%%%%%%%%%%%%%%%%%%%
%
One could try to implement the same idea for shifts in position
variables. Many problems arise however since the potential energy is
not a nicely regular function of the positions.

If one tries to consider shift distributions $\psi_N(\delta)$ that do
not depend on $Z_0$, then the limitation on $|\delta|$ is quite
drastic. In fact this example essentially works if only a fixed
(independent of $N$) number of coordinates $\delta_i$ are not $0$. 
   
Trying to generalize the second example by imposing that $Z_0+\delta$
and $Z_0$ are on the same energy surface also faces several problems. 
The main one is that the equation of the energy surface is not anymore
symmetric between the shift $\delta$ and the shift $-\delta$.

The only solution would be to write the equation only on the tangent
plane, {\it i.e.} something like
\[
\psi_N=\Psi(|\delta|)\,\delta_{\nabla
  H(Z_0)\cdot\delta=0} \, |\nabla H(Z_0)|.
\]
This is now nicely symmetric but poses other difficulties. For
instance one would need to make sure that $H(Z_0+\delta)\leq
H(Z_0)+C$. Expanding $H$, one would formally find a condition of
the type
\[
|\nabla^2 H(Z_0)|\,|\delta|^2\leq C.
\]
Unfortunately $\nabla^2 H(Z_0)$ is singular and in particular
unbounded. It is only bounded in average which would force us to
remove the initial conditions around which the energy has a singular
behavior.

Although this procedure could in principle be carried out
successfully, we do not wish to enter here into such technical
computations.  This essentially limits us to present explicit shift
distributions acting only on velocity variables.

However, let us point out that the evolution of the particles system
strongly mixes positions and velocities. Obviously, if we start from
two initial conditions $Z_0$ and $Z_0^\delta$ with same positions and
different velocities, we get at time $t>0$ configurations with
different positions.

So if we really need a shift distribution that acts also on the
positions at the origin $t=0$, a strategy may be to start at $t= -\tau
>0$, using a shift distribution acting only on velocities at this
time, and then to let evolves the particles till time $0$. The
original shift distribution transported by the flows is now a shift
distribution acting on position and speed. Using the Theorem
\eqref{theorem1}, we get that the average of these evolved shifts (in
a weak sense since we we are taking a logarithmic mean) is at most of
order $\delta_N$, provided the average of the original shifts was also
of order at most $\delta_N$.  Therefore by removing a set of vanishing
measure of initial conditions, one obtains a shift distribution in
positions and velocities that satisfies all the requirements.
%
%
%%%%%%%%%%%%%%%%%%%%%%%%%%%%%%%%%%%%%%%%%%%%%%%%
%
\section{Some useful bounds for the $6N$ dimensional $\mu_N$}
We shall make use of the following lemmas. Before stating them, we
introduce a notation for the projection of $\mu_N$ on the position
space. 
\begin{equation}
\nu_N(X^N) = \int \mu_N(X^N,V^N)\,dV^N \,.
\end{equation}
We also denote the $k$- marginal of $\nu_N$ by $\nu_N^k$.

%\begin{lemma}  $\exists A_N \subset [0,1]^{3N}$, 
%such that ${\cal L}^{3N}(A_N) \rightarrow 0$,
%  and for all $X_N \in (A_N)^c$,
%\begin{equation}
%E_{pot}(X_N) \geq (N-1)\,e_0 - CN^{(\alpha-1)/d}~.
%\end{equation}
%\label{lemma1}
%\end{lemma}

\begin{lemma}
For all $N$, we have:
\begin{equation}
 \left(\frac{2\pi}{\beta}\right)^{3N/2} \leq \mathcal{B}_N \leq
 \left(\frac{2\pi e^{- \beta \frac{\phi_{min}}3}}{\beta}\right)^{3N/2} . 
\end{equation}\label{lemma2}
\end{lemma}
We will also need the following estimate on the $k$th marginal of
$\mu_N$ defined by 
\begin{equation} \label{eq:marg}
\mu_N^k (Z^k) = \int \mu_N(dZ^{N-k})=
\frac{1}{\mathcal{B}_N}\;\int e^{-\beta H_{N}(Z^k,Z^{N-k})}\, dZ^{N-k} \,,
\end{equation}
and more precisely on the $k$ marginal $\nu_N^k$ of the projection
$\nu_N$ on the position variables.

\begin{lemma} \label{lemma3}
We define the constant $c_\beta = e^{-\beta \phi_{min}}$. Then, for
all $k$, we have   
\begin{eqnarray}
% \mu_N^k(Z^k) & \leq & c_\beta^k  \left(
% \frac\beta{2\pi}\right)^{\frac{3k}2} e^{E_{kin}(Z^k)}  \,, \\ 
\nu_N^k(X^k) & \leq & c_\beta^k
\end{eqnarray}
\end{lemma}

%Finally, we will use the fact that shifting $Z^N$ by a small amount 
%$\delta \in \mbox{supp} \, \psi_{N}$, changes the measure density at most by a
%multiplicative constant independent of $N$. More precisely

%\begin{lemma}\label{lemma4}
%For all $Z^N$ and all $\delta \in \mbox{supp}\,\psi_{N}$, with
%$\mbox{supp}\,\psi_{N}$ 
%satisfying \eqref{psi}, we have, with $K_{1}= e^{\beta \kappa}$ independent of
%$N$
%\begin{equation}
% K_{1}^{-1} \,d\mu_{N}(Z^N+\delta)\leq d\mu_{N}(Z^N)\leq K_{1}
%\,d\mu_{N}(Z^N+\delta)
%\end{equation}
%\end{lemma}

{\bf Proof of lemma \ref{lemma2}} To compute the integral defining
$\mathcal{B}_N$ \eqref{eq:betaN}, we may separate the integration in $X^N$
from the integration in $V^N$. In the $V^N$ variable, we have to
integrate a product of $3N$ independent real gaussians of variance
$\beta^{-1}$.  We obtain $(2\pi/\beta)^{3N/2}$. 

In the $X^N$ variable, we use Jensen inequality by the
convexity of exponential to get:
\begin{eqnarray*}
1
& = &  e^{\int \frac{-\beta}{2N}\sum_{i\neq
    j}\phi(|X^N_i-X^N_j|)\, d X^N  } \\ 
& \leq &  \int  e^{\frac{-\beta}{2N}\sum_{i\neq
    j}\phi(|X^N_i-X^N_j|)  }\,dX^N\\ 
& =& \left(\frac{2\pi}{\beta}\right)^{-3N/2}\;\int
e^{-\frac{\beta}{2}\sum_i |V_i^N|^2- \frac{\beta}{2N}\sum_{i\neq
    j}\phi(|X^N_i-X^N_j|) }\,dX^N\,dV^N\\
&=& \left(\frac{2\pi}{\beta}\right)^{-3N/2}\;\mathcal{B}_N,
\end{eqnarray*}
which gives the first inequality (We used that $\phi$ has zero
average). To obtain the second bound, 
it suffices to use the inequality  $E_{pot}(Z_0^N) \geq \frac{N}2 \phi_{min}$.
\cqfd

\medskip

{\bf Proof of lemma \ref{lemma3}} The proof follows the one introduced in 
\cite{MeSp} for the Lame-Enden equation.
As the measure $\mu_N$ factorizes in position and speed, we may write
\[
\nu_N(X^N) =  \frac{1}{\mathcal{B}_{N,X}} e^{-\beta  E_{pot}(X^N)} 
%, \quad \text{ with } 
%\beta_{N,x} = \int e^{-\beta E_{pot}(X^N)}\, dX^N
\]
Neglecting the terms in interaction energy involving (at least) one of
the first $k$ particles, we obtain 
\begin{eqnarray}
\nu_N^k (X^k) & = & \frac{1}{\mathcal{B}_{N,X}} \int e^{-\beta
  E_{pot}(X^k,X^{N-k})}\, dX^{N-k} \nonumber \\ 
& \leq & \frac{1}{\mathcal{B}_{N,X}}  e^{-\beta k  \phi_{min}}
\int e^{-\beta \frac{N-k}{N} E_{pot}(X^{N-k})} \, dX^{N-k}
~. \label{eq:bnu} 
\end{eqnarray}
the term $\frac{N-k} N$ being there because $E_{pot}(X_k) = (1/k) \sum_{i \neq
  j}^k \phi(|X_i-X_j|)$ for $k$ positions.  So we need an estimate on
terms of the kind  
\[
\Theta(k) = \int e^{-\beta \frac{k}{N} E_{pot}(X^k)} \, dX^k ~,
\]
We can relate this term to configurations with $k+1$ particles. 
First use Jensen inequality as the exponential is convex to get
\[
\begin{split}
\Theta(k) =& \int e^{-\beta \frac{k}{N}
  E_{pot}(X^k)}\,dX^k\\
=&\;\int \exp\left(-\int \left(\beta \frac{k}{N}
  E_{pot}(X^k)+\frac{\beta}{N} \sum_{i=1}^k
  \phi(|X_i-x_{k+1}|)\right)\;dx_{k+1}\right)\;
dX^k\\  
\leq & \int e^{-\beta \frac{k}{N} E_{pot}(X^k) -
    \frac{\beta}{N} \sum_{i=1}^k  \phi(|X_i-x_{k+1}|)} \, dX^k
  \,dx_{k+1} \\
 = & \int e^{-\beta \frac{k+1}{N} E_{pot}(X^{k+1})} \,
  dX^{k+1} \\ 
 \leq &  \Theta(k+1).
\end{split}\]
Since, $\Theta(N) =
\mathcal{B}_{N,X}$, we iterate this inequality and get $\Theta(N-k) \leq
 \mathcal{B}_{N,X}$. Putting this in \eqref{eq:bnu}, we get  
\[
\nu_N^k (X^k) \leq e^{ - k \beta \phi_{min}}  \,, 
\]
which is the result needed.
\cqfd

%\medskip
%{\bf Proof of lemma \ref{lemma4}}: the proof relies on the assumption
%that the perturbation $\delta$ to the configuration $Z^N$ does not
%change the energy 
%by more than an additive constant. Exponentiating \eqref{psi}, we have
%\[
%e^{-\beta \kappa}\,e^{-\beta H_{N}(Z_{N}+\delta)}\leq e^{-\beta H_{N}(Z_{N})} 
%\leq e^{\beta \kappa}\, e^{-\beta H_{N}(Z_{N}+\delta)}
%\]
%which immediately gives the result. \cqfd

%%%%%%%%%%%%%%%%%%%%%%%%%%%%%%%%%%%%%%%%%%%%%%%
\section{Proof of Theorem \ref{theorem1}}
%%%%%%%%%%%%%%%%%%%%%%%%%%%%%%%%%%%%%%%%%%%%%%%

During the course of the demonstration, $C$ will denote a constant (independent 
of $N$ and $\beta$), which value may change from line to line.

 From now on, we shall omit the superscript 
$N$ in the notation $Z^N$, as there will be no ambiguity. 
We have to estimate the derivative of $Q(t)$. Differentiating directly,
one obtains
\[\begin{split}
\frac{d}{dt} Q(t)\leq \int d\mu_N(Z_0)\int d\delta\, \psi_N(Z_0,\delta)\;
&\Bigg(\frac{\frac{1}{N}\sum_i |V_i-V_i^\delta|}{\delta_N
  +\|Z-Z^\delta\|_1}\\
&+\frac{\frac{1}{N^2}\sum_{i}
  |\sum_j (K(X_i-X_j)-K(X_i^\delta-X_j^\delta))|}{\delta_N
  +\|Z-Z^\delta\|_1}\Bigg),\\ 
\end{split}
\]
where $Z^\delta=(X^\delta,V^\delta)=Z(t,Z_0+\delta)$.

Note that the first term is obviously bounded by $1$ and hence
\[
\frac{d}{dt} Q(t)\leq 1 +\int d\mu_N(Z_0) \int d\delta\, \psi_N(Z_0,\delta)\;
\frac{\frac{1}{N^2}\sum_{i}
  |\sum_j (K(X_i-X_j)-K(X_i^\delta-X_j^\delta))|}{\delta_N
  +\|Z-Z^\delta\|_1}.
\] 
We define for a integer $L$ that will be fixed later
\begin{equation} \label{eq:Ci}
\mathcal{C}_i(Z_0,t) = \left\{ j\neq i~s.t.~|X_i(t)-X_j(t)|~\mbox{is among 
the}~L~\mbox{smallest}~|X_i-X_k| \right\}
\end{equation}
That is for each $i$, $\mathcal{C}_i$ collects the indices of particles
which are closest to particle~$i$ at time~$t$, following the flow. 
It also depends on the initial condition $Z_0$. We define similarly
$\mathcal{C}_i^\delta$. 

Accordingly, we decompose $dQ/dt$ as follows 
\[
\frac{d}{dt} Q(t)\leq C+S_1+S_1^\delta+S_2,
\]
with
\begin{equation}\begin{split}
&S_1=\int d\mu_N(Z_0) \int d\delta\, \psi_N(Z_0,\delta)\;
\frac{1}{\delta_N\, N^2}\sum_{i}\sum_{j\in \mathcal{C}_i\cup
  \mathcal{C}_i^\delta } 
|K(X_i-X_j)|,\\
&S_1^\delta=\int d\mu_N(Z_0) \int d\delta\, \psi_N(Z_0,\delta)\;
\frac{1}{\delta_N\, N^2}\sum_{i}\sum_{j\in \mathcal{C}_i\cup
  \mathcal{C}_i^\delta } 
|K(X_i^\delta-X_j^\delta)|,\\
&S_2=\int d\mu_N(Z_0) \int d\delta\, \psi_N(Z_0,\delta)\;
\frac{1}{N^2}\sum_{i}\sum_{j\not\in (\mathcal{C}_i\cup
  \mathcal{C}_i^\delta)}\frac{|K(X_i-X_j)-K(X_i^\delta-X_j^\delta)|}
{\delta_N+\frac{1}{N}\sum_i |X_i-X_i^\delta| }.  
\end{split}\end{equation}
%
%
%%%%%%%%%%%%%5555
\subsection{Bound on $S_1$}
%%%%%%%%%%%%%%%%%
Recalling the bounds on $K=-\nabla \phi$ from Th. \ref{theorem1}, one
simply begins with a discrete H\"older inequality for any $\gamma\leq 3$ 
\[\begin{split}
S_1\leq &\left(\int d\mu_N(Z_0) \int d\delta\, \psi_N(Z_0,\delta)\;
\frac{1}{\delta_N\, N^2}\sum_{i}\sum_{j\in \mathcal{C}_i\cup
  \mathcal{C}_i^\delta } 1 \right)^{1-\alpha/\gamma}\\
&\times\left(\int d\mu_N(Z_0) \int d\delta\, \psi_N(Z_0,\delta)\;
\frac{1}{\delta_N\, N^2}\sum_{i}\sum_{j\neq i}
\frac{1}{|X_i-X_j|^\gamma}\right)^{\alpha/\gamma}.   
\end{split}\]
We first use the fact that the integral of $\psi_{N}$ in $\delta$ is equal to
$1$ to get
\[
\int d\mu_N(Z_0) \int d\delta\, \psi_N(Z_0,\delta)\;
\frac{1}{\delta_N\, N^2}\sum_{i}\sum_{j\neq i}
\frac{1}{|X_i-X_j|^\gamma}=\int d\mu_N(Z_0)\;
\frac{1}{\delta_N\, N^2}\sum_{i}\sum_{j\neq i}
\frac{1}{|X_i-X_j|^\gamma}.
\]
We then perform the change of variable from $Z_{0}$ to $Z$, using the 
inverse flow of $Z(t,Z_0)$ which preserves the measure $\mu_N$; one finds
\[\begin{split}
\int d\mu_N(Z_0)\;
\frac{1}{\delta_N\, N^2}\sum_{i}\sum_{j\neq i}
\frac{1}{|X_i-X_j|^\gamma}\leq \frac{1}{\delta_N\,
  N^2}&\sum_{i}\sum_{j\neq i}\;\int d\mu_N(Z)\frac{1}{|X_i-X_j|^\gamma}.   
\end{split}\]
As the second marginal of $\mu_N$ is bounded by $c_\beta^2$ by Lemma
\ref{lemma3}, this implies
\[
\int d\mu_N(Z_0) \int d\delta\, \psi_N(Z_{0},\delta)\;
\frac{1}{\delta_N\, N^2}\sum_{i}\sum_{j\neq i}
\frac{1}{|X_i-X_j|^\gamma}\leq \frac{C\,c_\beta^2}{(3 -\gamma)\delta_N}.
\]
Hence
\[
S_1\leq \frac{C\,c_\beta^{2\alpha/\gamma}}{(3-\gamma) \delta_N}
\left(\int_{Z_0} d\mu_N(Z_0)\,\frac{1}{N^2}
\sum_i
(|\mathcal{C}_i|+|\mathcal{C}_i^\delta|)\right)^{1-\alpha/\gamma}.
\] 
To conclude, simply note that by definition
$|\mathcal{C}_i^\delta|=|\mathcal{C}_i|=L$ and so for any $a >
\frac{2\alpha}{3}$,  
\begin{equation}
S_1\leq \frac{C_a\,c_\beta^a}{\delta_N}\,
\left(\frac{L}{N}\right)^{1-\frac a 2} \,,
\end{equation} 
where $C_a$ satisfies $C_a \leq \frac{C}{3a -2 \alpha}$.
%%%%%%%%%%%%%%%%%%%%%%%%%
\subsection{Bound on $S_1^\delta$} \label{S1delta}
Using Fubini and the change of variable $Z_0 \mapsto Z_0 + \delta$,
and the image  
measure $\nu_N$, we may rewrite
\[
S_1^\delta =   \int \nu_N(Z_0) dZ_0 \,\;
\frac{1}{\delta_N\, N^2}\sum_{i}\sum_{j\in \mathcal{C}_i\cup
  \mathcal{C}_i^{-\delta} } |K(X_i-X_j)|
\]
And from the hypothesis \eqref{psi}, we may bound that that exactly
as the previous one. The only difference is that a constant $K_1$
appears and that we shall use $\beta'(N)$ instead of $\beta$. We get 
\[
S_1^\delta \leq \frac{C K_1 \,c_{\beta'(N)}^{2\alpha/\gamma}}{\delta_N}\,
\left(\frac{L}{N}\right)^{1-\alpha/\gamma} \,.
\]

%
%
%%%%%%%%%%%%%%%%%%%%%%%%%
\subsection{Bound on $S_2$}

By using the assumption on the second derivative of $\phi$ in
Th. \ref{theorem1}, one first bounds 
\[
|K(X_i-X_j)-K(X_i^\delta-X_j^\delta)|\leq
C\,(|X_i-X_i^\delta|+|X_j-X_j^\delta|)\;
\left(\frac{1}{|X_i-X_j|^{\alpha+1}}
+\frac{1}{|X_i^\delta-X_j^\delta|^{\alpha+1}}\right).
\]
Therefore defining the following matrix, with $\mathbb{I}_{A}$ the
characteristic function of the set $A$: 
\[
M_{ij}=\left(\frac{1}{|X_i-X_j|^{\alpha+1}}
+\frac{1}{|X_i^\delta-X_j^\delta|^{\alpha+1}}\right)\,({\mathbb{I}}_{j\not\in
  (\mathcal{C}_i\cup 
  \mathcal{C}_i^\delta)}+{\mathbb{I}}_{i\not\in
  (\mathcal{C}_j\cup 
  \mathcal{C}_j^\delta)}),
\]
one has
\[
S_2\leq \int d\mu_N(Z_0) \int d\delta\, \psi_N(Z_0,\delta)\;
\frac{\frac{1}{N^2}\sum_{i,j}M_{ij} |X_j-X_j^\delta|}
{\delta_N+\frac{1}{N}\sum_k |X_k-X_k^\delta| }.
\]
Consequently, if we use the classical matrix inequality $\|Mx\|_1 \leq \sup_{j}
\big( \sum_{i} |M_{ij}|) \big) \|x\|_1$, 
\[
S_2\leq \int d\mu_N(Z_0) \int d\delta\,
\psi_N(Z_0,\delta)\;\max_i\,\sum_j M_{ij}. 
\]
As before the terms in $M$ containing $|X_i^\delta-X_j^\delta|$ are the
equivalent of the ones with $|X_i-X_j|$, thanks to \eqref{psi}. 
Hence one has to bound
\[
S_2\leq C\,(S_2^1+S_2^2),
\]
with
\[\begin{split}
&S_2^1=\int d\mu_N(Z_0) \int d\delta\, \psi_N(Z_0,\delta)\;\max_i\,\sum_{j\not\in
  \mathcal{C}_i} \frac{1}{|X_i-X_j|^{\alpha+1}},\\
&S_2^2=\int d\mu_N(Z_0) \int d\delta\,
\psi_N(Z_0,\delta)\;\max_i\,\sum_{j\ s.t.\ i\not\in 
  \mathcal{C}_j} \frac{1}{|X_i-X_j|^{\alpha+1}}.\\
\end{split}\]
Since nothing depends on $\delta$ now, one may integrate
$\psi_N(Z_0,\delta)$ in $\delta$ with value $1$. 
Moreover changing variable from $Z_0$ to $Z$ (we
recall the flow is measure preserving), one simply finds
\[\begin{split}
&S_2^1=\int d\mu_N(Z) \;\max_i\,\sum_{j\not\in
  \mathcal{C}_i} \frac{1}{|X_i-X_j|^{\alpha+1}},\\
&S_2^2=\int d\mu_N(Z)\;\max_i\,\sum_{j\ s.t.\ i\not\in
  \mathcal{C}_j} \frac{1}{|X_i-X_j|^{\alpha+1}}.\\
\end{split}
\]
Let us now carefully bound each of these terms.

\medskip

{\bf The $S_2^1$ term}\\
We use 
\begin{equation}
\int f(X)~d\mu_N=\int_0^{+\infty}P(f(X)>l)~dl~,
\end{equation}
where $P$ is the probability with respect to the
measure $\mu_N$ on $\Pi^{3N}\times\R^{3N}$.  
We have to evaluate now expressions like
\begin{equation}
P \left( \max_i \frac{1}{N}\left( \sum_{j\notin\mathcal{C}_i} 
\frac{1}{|X_i-X_j|^{\alpha+1}}\right)>l \right)=
P\left( \exists i~s.t.~\frac{1}{N}\left( \sum_{j\notin\mathcal{C}_i} 
\frac{1}{|X_i-X_j|^{\alpha+1}}\right)>l\right)~.
\label{eq:P}
\end{equation}

To bound this probability, we will need the following lemma:
\begin{lemma} Given $x_1,\ldots,x_n\geq 0$ and $l\geq0$; given
  $(u_k)_{k=1}^n$ such that $\sum_{k=1}^n u_k=1$. If $\sum_{i=1}^n
  x_i>l$, then $\exists k\in[1,~n],~\exists
  i_1,\ldots,i_k~/~x_{i_r}>lu_k~,~\forall r=1,\ldots,k$.
\end{lemma}
{\emph{Proof:}} Let the $x_i$ be sorted $x_1\geq x_2\geq\ldots\geq x_n$, and
suppose the conclusion is not true. Then we have 
$$
x_1\leq l u_1,\ldots,x_n\leq l u_n~. 
$$
Thus $\sum_i x_i \leq l\sum_i u_k=l$.\cqfd

\medskip

We apply the lemma to Eq.~\ref{eq:P}, with $u_k=\frac{C_N \nu}k
\left(\frac k N \right)^\nu$ and 
$0<\nu<1$ to be determined. $C_N$ is chosen such that
\mbox{$\sum_{k=1}^N u_k=1$.} Using Riemann sums, we see that
$\lim_{+\infty} C_N = 1$, and that $C_N \geq 1, \forall N \geq 1$. Hence we get 
\begin{eqnarray*}
P_l & = & P \left( \exists i~s.t.~\frac{1}{N}  \left( \sum_{j\notin\mathcal{C}_i} 
      \frac{1}{|X_i-X_j|^{\alpha+1}}\right) >l \right) \\
  &  \leq  & \sum_{k=1}^{N-L} P\left( \exists i;j_1,\ldots,j_k \notin
      \mathcal{C}_i,~
      \frac{1}{N}\frac{1}{|X_i-X_{j_r}|^{\alpha+1}}>l\frac{C_N \nu}{k}(k/N)^\nu
      \right) \\ 
 &  \leq & \sum_{k=1}^{N-L} P\left( \exists i;j_1,\ldots,j_k \notin
      \mathcal{C}_i,~ |X_i-X_{j_r}|<\left(\frac1{ \nu l
      }\right)^{1/(\alpha+1)} \left(\frac k N \right)^{\lambda/3} 
      \right) \\ 
 &  = & \sum_{k=1}^{N-L}P_{l,k} ~,
\label{eq:Pmu}
\end{eqnarray*}
where for simplicity we have introduced the parameter $\lambda=3
(1-\nu)/(\alpha+1)$.  

To estimate the probability $P_{l,k}$, once the particle $i$ is
chosen, we have a constraint on the position of $k$ particles, 
which have to be close enough to particle~i, plus constraints on the
position of $L$ distinct particles, from the definition of
$\mathcal{C}_i$. This event concern $k+L+1$ particles, and to estimate
it, we will use an estimate of its volume $P_{l,k}^u$ in the
configuration space $\T^{3(k+L+1)}\times\R^{3(k+L+1)}$. It thus
involves the $(k+L+1)$ marginal of $\mu_N$ which is bounded by
$c_\beta^{k+L+1}$ by Lemma \ref{lemma3}.
 
This leads to the following estimates
\[
P_{l,k} \leq  C\;c_\beta^{k+L+1}\,   N~
\left(\frac{1}{\nu l}\right)^{\frac{3(k+L)}{\alpha+1}} \left(
\frac{k}{N}\right)^{\lambda(k+L)}.\\ 
\]  
Moreover, using a simplified version of Binet  formula (See \cite{WW})
\[
n! = \sqrt{2\pi} n^{n+ \frac12} e^{-n + \frac{\theta}{12n}} \,,
\text{ for some } \theta \in (0,1) 
\]
the binomial coefficient $C_N^p$ may be bound by: 
\begin{eqnarray}
C_N^p  & = & \frac{N!}{p!(N-p)!} =  (2\pi)^{-\frac12}
e^{\frac{\theta_N}{12N} - \frac{\theta_p}{12p} -
  \frac{\theta_{N-p}}{12(N-p)}} \frac{N^N}{p^p (N-p)^{N-p}} 
\left(\frac{N}{p(N-p)}\right)^{1/2} \\ 
& \leq &  (2\pi)^{-\frac12} e^{\frac1{12N}}
\left(\frac{N}{N-1}\right)^{1/2} \left( \frac{N}{p}\right)^p \left(1
+ 
\frac{p}{N-p}\right)^{N-p} \leq  \left( \frac{N e}{p}\right)^p ~. 
\end{eqnarray}
And we do not forget that since $C_N^p = C_N^{N-p}$ we may use the
same inequality with $p$ replaced by $N-p$. 
Inserting this in the above inequality, we get:
\begin{eqnarray*}
P_{l,k} & \leq & C c_\beta^{k+L+1} N \left( \frac{N}{k+L}\right)^{k+L}
e^{k+L} \left(\frac{1}{\nu l}\right)^{\frac{3(k+L)}{\alpha+1}}\left(
\frac{k}{N}\right)^{\lambda(k+L)} \\ 
& \leq & C c_\beta \, N \,\left( \frac{k}{N}\right)^{(\lambda-1)(k+L)}
\left(\frac{A c_\beta'}{\nu l}\right)^{\frac{3(k+L)}{\alpha+1}} \\ 
& \leq &  C c_\beta \, \left( \frac{k}{N}\right)^{(\lambda-1)(k+L) - 1}
k  \, \left(\frac{A c_\beta'}{\nu l}\right)^{\frac{3(k+L)}{\alpha+1}}~,
\\ 
\end{eqnarray*}
where $c_\beta' = c_\beta^{(\alpha+1)/3}$, and $A= e^{(\alpha+1)/3}$
is a numerical constant.  Now taking $\nu$ close enough to $0$
(precisely $\nu < \frac{2- \alpha}3$), one 
has $\lambda
>1$ 
and then we take as well $L \geq (\lambda -1)^{-1}$ (recall that $L$
has yet to be fixed)~; hence
\[
P_{l,k}\leq  C c_\beta \, k \, \left(\frac{A
  c_\beta'}{\nu l}\right)^{\frac{3(k+L)}{\alpha+1}}.
\]  
If we sum on $k$, we get:
\begin{eqnarray}
\sum_{k=0}^{N-L} P_{l,k} & \leq & C c_\beta \, \left(\frac{A
  c_\beta'}{ \nu l}\right)^{\frac{3L}{\alpha+1}} \, \sum_{k=0}^{N-L} k
\,\left(\frac{A c_\beta'}{\nu l}\right)^{\frac{3k}{\alpha+1}} \\ 
& \leq &  C \, c_\beta \, \left(\frac{A
  c_\beta'}{\nu l}\right)^{\frac{3(L+1)}{\alpha+1}} \frac{1}{(1- (A
  c_\beta'/ \nu l)^{3/(\alpha+1)})^2}  \leq C \,  
\, c_\beta\, \left(\frac{A c_\beta'}{ \nu l}\right)^{\frac{3(L+1)}{\alpha+1}}~,
\end{eqnarray}
provided $l \geq l_0=\frac{2  \, A \, c_\beta'}{ \nu}$. If we take moreover 
$L \geq p \geq p(\alpha+1)/3$ for some $p \geq 0$, we get a simpler bound: 
\[
\sum_{k=0}^{N-L} P_{l,k}  \leq C \, c_\beta \left(\frac{A \,
  c_\beta'}{ \nu l}\right)^p ~. 
\]
Remark that the conditions on $L$ depend only on $\alpha$ and
$\lambda$ (which depends only on $\alpha$). In particular, those
conditions are 
independent of the parameter $\beta$. Thus, we have 
\begin{eqnarray}
P\left( \exists i~s.t.~\frac{1}{N}\left( \sum_{j\notin\mathcal{C}_i} 
      \frac{1}{|X_i-X_j|^{\alpha+1}}\right)>l\right) & \leq & 1
      ~\text{for}~l<l_0 = \frac{2A}{ \nu} c_\beta^{\frac{\alpha+1}3}    \\ 
& \leq & C c_\beta^{1+ p \frac{\alpha+1}3}( \nu l)^{-p} ~\text{for}~l\geq l_0' 
\end{eqnarray}

Integrating this quantity in $l$, one obtains for  $M$ such that  $C
c_\beta^{1+ p \frac{\alpha+1}3}( \nu M)^{-p}  =  1$ that 
\begin{eqnarray*}
\int_0^\infty P(\dots>l) \,dl & =  & \int_0^M P(\dots>l) \,dl +
\int_M^\infty P(\dots>l) \,dl \\ 
& \leq & M + C c_\beta \left(\frac{A c_\beta'}{ \nu}\right)^p
\int_M^\infty  \frac{dl}{l^p} \\ 
& \leq & M + C c_\beta \left(\frac{A c_\beta'}{ \nu}\right)^p
\frac{M^{1-p}}{p-1} \\ 
& = & \frac{p}{p-1} M \leq \frac{C}\nu  ~c_\beta^{\frac1p+  \frac{\alpha+1}3}
%\int_0^\infty P(\dots>l) \,dl & \leq & C c_\beta^\kappa + C
%c_\beta^{1- (p-1)\kappa} c_\beta'^{p} ~,   \\ 
%& \leq & C c_\beta^{\frac{\alpha+1}{3} + \frac{1}{p}} \leq C c_\beta ~,
\end{eqnarray*}
for $p \geq \frac{3}{2-\alpha}$. And finally, we get that
\[
S_{2}^1 \leq \frac{C}\nu ~ c_\beta ~,
\]
for some constant $C$, provided that 
$L$ is large enough.

Since every calculation have been performed, we see that a possible
for $\nu$ is  
$\nu = \frac{2-\alpha}{6}$ in which case the condition on $L$ is exactly
\[
L \geq \frac{6}{2-\alpha} 
\]
With that choice of $\nu$, we get
\[
S_{2}^1 \leq \frac{C}{2 -\alpha} ~ c_\beta ~,
\]

\medskip

{\bf The $S_2^2$ term}:
Through the same type of computations, we are led to evaluate expressions like

\begin{eqnarray}
P_l' & = & P\left( \exists j~s.t.~\frac{1}{N}\left( \sum_{i/j\notin\mathcal{C}_i} 
      \frac{1}{|X_i-X_j|^{\alpha+1}}\right)>l\right) \\
  & \leq & \sum_{k=1}^N P\left( \exists j;i_1,\ldots,i_k~/~j\notin
      \mathcal{C}_{i_r}
      \frac{1}{N}\frac{1}{|X_j-X_{i_r}|^{\alpha+1}}>l\frac{C_N
        \nu}{k}\left(\frac k N\right)^\nu\right) 
      \\
  & \leq & \sum_{k=1}^N P\left( \exists j;i_1,\ldots,i_k~/~j \notin
      \mathcal{C}_{i_r} |X_j-X_{i_r}|< \frac1{(\nu l)^{\alpha+1}}
    \left(\frac k N \right)^{(1-\nu)/(\alpha+1)}\right) \\
    & = & \sum_{k=1}^N P'_{l,k} ~.
\label{eq:Pmu2}
\end{eqnarray}
Since the sum is performed on the particles $i$ such that $j \notin
  \mathcal{C}_i$, we cannot choose $L+k$ particles close to $j$ as for
  $S_2^1$. But,
  we have nevertheless to choose $k$ particles close enough to $j$, a
  probability that will give a good bound if $k \geq L$. If $k \leq
  L$, once a particle $i$ close to $j$ such that $j \notin
  \mathcal{C}_i$ is chosen, one knows that there exist $L$ other
  particles close to $j$.
This will be enough to bound the probability. 
   
  In the second case ($k \leq L$), we pick up a particle $j$ ($N$
  possibilities), at least another particle $i$ (since $k\geq 1$)
and then we have to choose $L$ other particles closer to
  $i$ than $j$ is. Since $|X_j-X_i|$ has to be less than
$C\, l^{-1/(\alpha+1)}\,(k/N)^{(1-\nu)/(\alpha+1)}$, 
\begin{eqnarray}
  \sum_{k=0}^L P'_{k,l} & \leq & \sum_{k=0}^L c_\beta^{L+2} N^2 C_N^L
  \left(\frac{1}{\nu l}\right)^{3L/(\alpha+1)}  
\left(\frac{k}{N}\right)^{\lambda L}\\
& \leq & C  \sum_{k=0}^L c_\beta^2 \, N^2
\,\left(\frac{L}{N}\right)^{(\lambda-1)L} \left(\frac{A
  c_\beta'}{\nu l}\right)^{3L/(\alpha+1)} \\ 
& \leq &  C c_\beta^2  \, N^3 \,
\left(\frac{L}{N}\right)^{(\lambda-1)L} \left(\frac{A
  c_\beta'}{\nu l}\right)^{3L/(\alpha+1)}~. 
\end{eqnarray}
If $L \leq \sqrt N$, we may use
\[
N ^3 \left( \frac{L}{N} \right)^{(\lambda-1)L} \leq  N^{3 -
  \frac{\lambda-1}2 L} \leq 1 \,, 
\]
as soon as $L \geq \frac6{\lambda-1}$. In that case, 
we may use this bound in last inequality and obtain:
\begin{eqnarray*}
\sum_{k=0}^L P'_{k,l} & \leq & C c_\beta^2  \, \left(\frac{A
  c_\beta'}{l}\right)^{3L/(\alpha+1)} \\ 
& \leq & C c_\beta^2 \left(\frac{A c_\beta'}{l}\right)^{3L/(\alpha+1)} ~,
\end{eqnarray*}
if $L \leq \sqrt N $.

In the first case $k>L$, we pick up the particle $j$, and then
choose $k$ particles $i_r$ close to $j$. We obtain as previously
\begin{eqnarray*}
\sum_{k=L+1}^{N-1} P'_{k,l} & \leq & \sum_{k=L+1}^{N} c_\beta^{k+1} N C_N^k
\left(\frac{1}{\nu l}\right)^{3k/(\alpha+1)}  
\left(\frac{k}{N}\right)^{\lambda k}\\
& \leq & C  \sum_{k=L+1}^{N} c_\beta \, N
\,\left(\frac{k}{N}\right)^{(\lambda-1)k} \left(\frac{A
  c_\beta'}{\nu l}\right)^{3k/(\alpha+1)} \\ 
& = &  C \, c_\beta  \sum_{k=L+1}^{N}    k \,
\left(\frac{k}{N}\right)^{(\lambda-1)k - 1} \left(\frac{A
  c_\beta'}{l}\right)^{3k/(\alpha+1)} \\ 
& \leq & C c_\beta \sum_{k=L+1}^{N}  k\left(\frac{A
  c_\beta'}{\nu l}\right)^{3k/(\alpha+1)}  \\ 
& \leq & C c_\beta \left(\frac{A c_\beta'}{\nu l}\right)^{3L/(\alpha+1)} ~,
\end{eqnarray*}
where we again restricted ourselves to $(\lambda-1)L \geq 1$ and
assume $l \geq l_0$ .  
% Note again that as before since $l$ is large the dominant term is the
% geometric sum with $(A\,c_\beta'/l)^{3k/(\alpha+1)}$.
Putting the two sum together, we get the bound
\[
\sum_{k=1}^{N-1} P'_{k,l}  \leq   C  c_\beta^2 \left(\frac{A
  c_\beta'}{\nu l}\right)^{3L/(\alpha+1)} 
\]

It remains to integrate in $l$.  Doing exactly as for the $S_2^1$
term, and choosing the same $\nu$, we will get  
\[
S_2^2\leq \frac C {2 -\alpha}\,c_\beta.
\]
The only difference is that it will require $p \geq
\frac{6}{2-\alpha}$ and thus $L \geq \frac{12}{2-\alpha}$ 
%
%%%%%%%%%%%%%%%%%%%%%%%%5
\subsection{Conclusion of the proof}
%%%%%%%%%%%%%%%%%%%%%%%%
Putting all together, we finally we may bound
\[
\frac{dQ}{dt} 
\leq 1 +C_a\,\frac{c_\beta^a}{\delta_N}\;
\left(\frac{L}{N}\right)^{1-\frac a 2 }+ \frac{C}{2-\alpha}\,c_\beta, 
\]
with $c_\beta = e^{- \beta \phi_{min}}$ and $a> \frac{2\alpha}3$,
where $L$ is subject to the restrictions 
(with  the choice $\nu= \frac{2-\alpha}6$ which means that $\lambda =
\frac{4+\alpha}{2+2\alpha}$ ) 
\[
L\geq \frac{36}{2-\alpha}, \quad ,\quad L\leq \sqrt N.
\]

It is possible only if $N \geq \frac{6^4}{(2-\alpha)^2}$ and in that
case it is clear that $L$ should be chosen as small as possible and
from 
the constraint that means
\[
L=\frac{36}{2- \alpha}.
\]
 With this choice, one has
\[
\frac{dQ}{dt} 
\leq 1+\frac{C_a}{2-\alpha}\,c_\beta^a\;
\frac{1}{\delta_N N^{1-\frac a 2 }}+\frac C{2-\alpha}\,c_\beta.
\]
Now if one takes $\delta_N = N^{-\varepsilon}$, we can get
a uniform bound in $N$, only if
\[
\varepsilon \leq 1-\frac a 2 .
\]
If this is true, we get 
\[
\frac{dQ}{dt} 
\leq 1+ \frac{C_{a}\,c_\beta^{a}\,+C\,c_\beta}{2-\alpha}.
\]
with $C_a \leq \frac{C}{3a -2\alpha}$ which is the result given by
Theorem \ref{theorem1}. 
%%%%%%%%%%%%%%%%%%%%%%%%%%%%%%%%%
%%%%%%%%%%%%%%%%%%%%%%%%%%%%%%%%%%%%
\medskip

{\bf Acknowledgments:} this work was partially supported by the ANR 
09-JCJC-009401 INTERLOP project.


\begin{thebibliography}{10}
%%%%%%%%%%%%%%%%%%%%%%%%%%%%%%%%%%
%%%%%%%%%%%%%%%%%%%%%%%%%%%%%%%%%%
\bibitem{Am} L. Ambrosio, {\em Transport equation and Cauchy problem
    for $BV$ vector fields}.  Invent. Math.  {\bf 158},  no. 2, pp
  227--260, 2004.

\bibitem{ALM} L. Ambrosio, M. Lecumberry, S. Maniglia,
 Lipschitz regularity and approximate differentiability of the
 DiPerna-Lions flow. {\em Rend. Sem. Mat. Univ. Padova} {\bf 114} (2005),
 29--50.

%\bibitem{BGM00} C. Bardos, F. Golse, N. Mauser \emph{Weak coupling
%    limit of the N-particle Schr\"odinger equation}, Methods and
%  Applications of Analysis {\bf 7}, pp 275-293, 2000.

%\bibitem{BGMY} C. Bardos, L. Erd\"os, F. Golse, N. Mauser, H.T. Yau,
%  {\em Derivation of the Schr�dinger-Poisson equation from the quantum
%    $N$-body problem}. 
%  C. R. Math. Acad. Sci. Paris  {\bf 334},  no. 6, pp 515--520, 2002. 

\bibitem{Bo} F. Bouchut, Renormalized solutions to the Vlasov 
equation with coefficients of bounded variation. {\em
 Arch. Ration. Mech. Anal.} {\bf 157} (2001), pp. 75--90.

%\bibitem{BC} F. Bouchut, G. Crippa, Uniqueness, renormalization, 
%and smooth approximations for linear transport equations.  {\em SIAM
%J. Math. Anal.}  {\bf 38}  (2006),  no. 4, 1316--1328. 

\bibitem{BH77} W. Braun and K. Hepp, \emph{The Vlasov dynamics and its 
fluctuations in the 1/N limit of interacting particles}, Comm. Math. Phys.
{\bf 56}, pp 101-113, 1977.

\bibitem{CIP} C. Cercignani, R. Illner and M. Pulvirenti, {The
mathematical theory of dilute gases}, {Applied Mathematical Sciences},
106, {\em Springer-Verlag New-York},
1994.

\bibitem{CJ} N. Champagnat, P.E. Jabin, { \em 
Well posedness in any dimension for some hamiltonian flows with
  non $BV$ force terms}. To appear in
  {Comm. Partial Differential Equations}.

\bibitem{CD} G. Crippa, C. DeLellis, Estimates and regularity results
  for the DiPerna-Lions flow. 
{\em J. Reine Angew. Math.} {\bf 616} (2008), 15--46.

\bibitem{DeL} C. De Lellis, Notes on hyperbolic systems of 
conservation laws and transport equations. Handbook of differential
equations, Evolutionary equations, Vol. 3 (2007).

\bibitem{DL} R.J. DiPerna, P.L. Lions, Ordinary differential
  equations,  transport theory and Sobolev spaces. {\em Invent. Math.} {\bf 98}
(1989), 511--547. 

\bibitem{Dobrushin} R. L. Dobrushin, {\em Vlasov
    equations}, Funct.  Anal. Appl. {\bf 13}, pp 115-123, 1979.

\bibitem{GHL90} J.~Goodman, T.~Y.~Hou and J. Lowengrub, \emph{Convergence of 
the point vortex method for the 2-D Euler equations}, Comm. Pure Appl. Math.
{\bf 43}, pp 415-430, 1990.

\bibitem{Hau}  M. Hauray, {\em On Liouville transport equation with
    force field in $BV\sb {\rm loc}$},  
Comm. Partial Differential Equations  {\bf 29},  no. 1-2, pp 207--217, 2004.

\bibitem{Hau2} M. Hauray, {\em 
Wasserstein distances for vortices approximation of
Euler-type equations}.  Math. Models Methods Appl. Sci.  {\bf 19}  (2009),
no. 8, 1357--1384.

\bibitem{HJ06} M. Hauray and P.~E.Jabin, \emph{N-particles approximation of 
the Vlasov equation with singular potential},
Arch. Ration. Mech. Anal.  {\bf 183},  no. 3, pp 489--524, 2007.

\bibitem{HLL} M. Hauray, C. Le Bris, P.L. Lions, Deux remarques 
sur les flots g\'en\'eralis\'es d'\'equations diff\'erentielles
ordinaires. {\em C. R. Math. Acad. Sci. Paris}  {\bf 344}  (2007),
no. 12,
 759--764.

%\bibitem{IP} R. Illner and M. Pulvirenti, {\em Global validity of
%the Boltzmann equation for two and three-dimensional rare gas in
%vacuum}, Comm. Math. Phys. {\bf 121}, pp 143--146, 1989.

\bibitem{Ja} P.E. Jabin, {\em Differential Equations with singular
  fields}. Preprint.

\bibitem{JO} P.E. Jabin, F. Otto, {\em Identification of the dilute
regime in particle sedimentation}, Comm. Math. Phys., {\bf 250}, 
pp 415--432, 2004.

%\bibitem{lanford1} O.E. Lanford,  {\em 
%Dynamical systems, theory and applications} (Battelle
%Rencontres, Seattle, Wash., 1974), pp. 1--111, Lecture Notes in Phys.,
%Vol. 38, Springer, Berlin, 1975.

\bibitem{Lanford} O.E. Lanford, {\em On a derivation of the Boltzmann
    equation}.  International Conference on Dynamical Systems in
  Mathematical Physics (Rennes, 1975),  pp. 117--137. Asterisque,
  No. 40, Soc. Math. France, Paris, 1976.  

\bibitem{MeSp} J. Messer, H. Spohn, {\em 
Statistical mechanics of the isothermal Lane-Emden equation}.
J. Statist. Phys. 29 (1982), no. 3, 561--578. 


\bibitem{Neunzert} H. Neunzert, J. Wick, \hbox{Theoretische und 
numerische Ergebnisse zur nicht} linearen Vlasov Gleichung. 
{\em \hbox{Numerische L\"osung nichtlinearer partieller Diffe-} rential und
Integrodifferentialgleichungen (Tagung, Math. Forschungsinst.,
Oberwolfach, 1971)}, pp. 159--185. Lecture Notes in Math., Vol. 267, 
Springer, Berlin, 1972. 

\bibitem{Sch1} S. Schochet, {\em The weak vorticity formulation
of the 2-D Euler equations and concentration-cancellation}, Comm.
Partial Differential Equations, {\bf 20}, pp 1077-1104, 1995.

\bibitem{Sch2} S. Schochet, {\em The point-vortex method for
periodic weak solutions of the 2-D Euler equations},  Comm. Pure
Appl. Math., {\bf 49}, pp 911-965, 1996.

\bibitem{Sp} H. Spohn, {Large scale dynamics of interacting
particles}, {\em Springer-Verlag Berlin}, 1991.
\bibitem{VA} H.D. Victory, jr., E.J. Allen, {\em The convergence
theory of particle-in-cell methods for multidimensional Vlasov-Poisson
systems}, SIAM J. Numer. Anal., {\bf 28}, pp 1207--1241, 1991.

\bibitem{WW}   E. T. Whittaker and  G.N. Watson, {\em A course of
  modern analysis}, Cambridge University Press, 1927. 
   

\bibitem{Wo} S. Wollman, {\em On the approximation of the
Vlasov-Poisson system by particles methods}, SIAM J. Numer. Anal.,
{\bf 37},
pp 1369--1398, 2000.

%\bibitem{eqstatmech} Je pense ici \`a des articles sur la m\'eca stat
%d'\'equilibre  
%de syst\`emes de particules avec interactions singuli\`eres (Messer-Spohn 82?, 
%Caglioti -Lions-Marchioro-Pulvirenti, Kiessling-Lebowitz... )
%\bibitem{Kiessling07} papier annon\c{c}ant le r\'esultat, sans suite...
\end{thebibliography}
\end{document}